\documentclass[letterpaper, 10 pt, conference]{ieeeconf}  % Comment this line out
\IEEEoverridecommandlockouts                              % This command is only
\usepackage{graphics} % for pdf, bitmapped graphics files
\usepackage{epsfig} % for postscript graphics files
\usepackage{amsmath} % assumes amsmath package installed
\usepackage{amssymb}  % assumes amsmath package installed
\usepackage{color}
\usepackage{subcaption}

\usepackage{amsthm}
\usepackage{caption}
\usepackage{epstopdf}
\let\labelindent\relax
\usepackage{enumitem}

\newlist{Aenumerate}{enumerate}{1}
\setlist[Aenumerate]{label=B.\arabic*}

\newcommand{\mcl}[1]{\mathcal{ #1}}
\newcommand{\mbf}[1]{\mathbf{ #1}}

\newtheorem{thm}{Theorem}
\newtheorem{defn}{Definition}[thm]
\newtheorem{lem}{Lemma}[thm]

\newcommand{\mylen}{3.1pt}
\newcommand{\hinfty}{\ensuremath{H_{\infty}}}
\newcommand{\ip}[2]{\langle{#1},{#2}\rangle}
\newcommand{\vmatwo}[2]{\begin{bmatrix}
		#1 \\
		#2
\end{bmatrix}}

\renewcommand{\th}{\ensuremath{\theta}}
\newcommand{\et}{\eta}

\newcommand{\mc}[1]{\ensuremath{\mathcal{#1}}}

\newcommand{\opPb}{\ensuremath{\mc{P}_{\small \{P,Q,R_1,R_2\}}}} 
\newcommand{\opPa}{\ensuremath{\mc{P}}}
\newcommand{\opPmn}{\ensuremath{\mc{P}_{\{M,N_1,N_2\}}}}
\newcommand{\myint}{\int_{0}^{L}}
\newcommand{\myinta}[1]{\int_{0}^{#1}}
\newcommand{\myintb}[1]{\int_{#1}^{L}}

\newcommand{\ltwo}{\ensuremath{L_2}}
\newcommand{\rl}{\ensuremath{\mathbb{R}}}

\newcommand{\norm}[1]{\ensuremath{\Vert #1\Vert}}
\newcommand{\sob}[1]{\ensuremath{W^{2,n}(#1)}}

\newcommand{\bmat}[1]{\begin{bmatrix}#1\end{bmatrix}}

\newcommand{\R}{\mathbb{R}}

\allowdisplaybreaks

\title{\LARGE \bf
Computing Input-Ouput Properties of Coupled Linear PDE systems
}
\author{%
	Sachin Shivakumar$^{1}$\\
	\and 
	Matthew M. Peet$^{1}$
	\thanks{$^{1}$ Sachin Shivakumar \{sshivak8@asu.edu\} and Matthew M. Peet \{mpeet@asu.edu\} are with School for Engineering of Matter, Transport and Energy, Arizona State University, Tempe, AZ, 85298 USA}
}
\begin{document}
\maketitle

\setlength{\abovedisplayskip}{\mylen}
\setlength{\belowdisplayskip}{\mylen}
%%%%%%%%%%%%%%%%%%%%%%%%%%%%%%%%%%%%%%%%%%%%%%%%%%%%%%%%%%%%%%%%%%%%%%%%%%%%%%%%
\begin{abstract}
In this paper, we propose an LMI-based approach to analyze input-output properties of coupled linear PDE systems. This work expands on a newly developed state-space theory for coupled PDEs and extends the positive-real and bounded-real lemmas to infinite dimensional systems. We show that conditions for passivity and bounded \ltwo ~gain can be expressed as linear operator inequalities on $\rl\times\ltwo$. A method to convert these operator inequalities to LMIs by using parameterization of the operator variables is proposed. This method does not rely on discretization and as such, the properties obtained are prima facie provable. We use numerical examples to demonstrate that the bounds obtained are not conservative in any significant sense and that the bounds are computable on desktop computers for systems consisting of up to 20 coupled PDEs. 
\end{abstract}

%%%%%%%%%%%%%%%%%%%%%%%%%%%%%%%%%%%%%%%%%%%%%%%%%%%%%%%%%%%%%%%%%%%%%%%%%%%%%%%%
\section{INTRODUCTION}\label{sec:intro}
%Many of engineering systems have states that depend on multiple variables. Dynamics of such systems generally involve Partial Differential Equations (PDEs). Such systems,
Partial Differential Equations (PDE) are used to model systems whose state varies not just in time, but also depend on one or more independent variables. For example, PDEs are used to model systems that have deformable structures \cite{Krstic2013}, thermo-fluidic interactions \cite{henningson2009}, and chemical processes \cite{PANJA2012, christofides96}.
Furthermore, the states of these PDEs are often vector-valued, representing, e.g. changes in temperature due to flow or interaction between chemical subspecies.

In this paper, we seek to develop algorithms which establish provable properties of linear, coupled PDE systems with inputs and outputs. Specifically, we develop Linear Matrix Inequality (LMI) tests for passivity and $L_2$-gain of PDE systems where for $w \in L_2^m[0,\infty)$, the system is defined by solutions to the following set of equations
\setlength{\arraycolsep}{2pt}
\begin{align}
  &\dot{x}(s,t) = A_0(s)x(s,t) + A_1(s)x_s(s,t)\notag \\
  &\qquad\qquad\qquad+A_2(s)x_{ss}(s,t) + B_1(s)w(t),\notag \\
  &y(t)= C_1 z(t) + \int_{0}^L \left(C_a(s)x(s,t)+C_b(s)x_s(s,t)\right)ds,\notag \\
  &B\bmat{x(0,t)& x(L,t)& x_s(0,t)& x_s(L,t)}^T = 0, ~ x(s,0)=0, \label{eq:PDE0}
 \end{align}
where if $B\in \R^{2n \times 4n}$ has row rank of $2n$, then $x(s,t) \in \R^n$ and $y(t)\in \R^q$ are uniquely defined. This system has a distributed input (typically modeling disturbances), and a combined boundary-valued/distributed output. Our goal, then, is
\begin{enumerate}
\item \textbf{$L_2$ Gain:} To find the smallest $\gamma$ such that $\norm{y}_{L_2}\le \gamma \norm{w}_{L_2}$ for all $w \in L_2^m$ and \item \textbf{Passivity:} To check wehther $\ip{y}{w}_{L_2}\ge 0$ for all $w \in L_2^q$.
\end{enumerate}

Most methods for analysis of PDE systems involve approximating the continuous infinite-dimensional state variables by a finite set of states \cite{prud2002reliable}, \cite{ELFARRA2003715} - yielding a system defined by ODEs. These methods, although well-studied, are limited by the fact that properties proven for the ODE approximation of a PDE are not prima facie provable for the PDE - although in some cases a \textit{posteriori} error bounding may be used to obtain properties such as $L_2$ gain bounds for the original PDE. Furthermore, a \textit{posteriori} error bounding will typically depend on the method and level of discretization and may involve substantial conservatism.Besides, ODE approximations of PDE models often require large number of states, resulting in intractably large optimization problems when analyzed in the LMI framework.

%In \cite{krstic2008boundary}, a backstepping method was used for analysis and control of PDEs.
Some prior work on properties of PDEs in an infinite-dimensional framework includes~\cite{GAYE2013} and~\cite{fridman} which proposed LMIs for \hinfty ~analysis of parabolic and hyperbolic PDEs, but were restricted to PDE systems with a single state, i.e. $n=1$. Other works, such as ~\cite{ahmadi2014} and~\cite{ahmadi16dissi}, proposed LMIs for \ltwo ~gain and passivity analysis of PDE systems that resulted in less conservative bounds, but %the bounds therein were still conservative and, moreover, 
even for small-scale linear problems, the resulting LMIs were significantly larger than the LMIs we use. Also note, the methods mentioned here restrict to PDEs with one spatial dimension, but there are other methods that use LMIs for analysis of PDEs with N spatial dimensions, such as \cite{fridman2d}.

Our approach is based on a generalization of the LMI framework for analysis of ODE systems to infinite-dimensional systems. Specifically, the LMI framework uses positive matrix variables to parameterize quadratic Lyapunov functions for analysis and control of ODE systems. In our approach, we use linear operators parameterized by matrix-valued polynomials to parameterize quadratic Lyapunov functionals for infinite-dimensional systems. This is an extension of the work on stability analysis in~\cite{peet2018}. Note that such an approach was previously used for time-delay systems (e.g. in \cite{peetcdc2017}), but has not been extended to PDEs. 

Here, we briefly recall the LMI approach to bound the \hinfty ~norm of an ODE system. For an ODE system  represented in traditional state-space representation \eqref{eq:ODE},
\begin{align}\label{eq:ODE}
\dot{x}(t) &= Ax(t) + Bw(t),~x(0)=0\nonumber\\
y(t) &= Cx(t) + Dw(t)
\end{align}
the following LMI condition \cite{boyd1994linear}, established using bounded-real lemma, can be used to find a bound on \hinfty ~norm.
\begin{thm}
Define: \begin{align*}
    G(s) := C(sI-A)^{-1}B+D.
\end{align*}
If there exists a positive definite matrix $P$, such that
\begin{align}\label{eq:lmihinf}
\begin{bmatrix}
A^TP+PA & PB & C^T\\
B^TP & -\gamma I & D^T\\
C & D & -\gamma I
\end{bmatrix}\le 0,
\end{align}
then $\norm{G}_{\hinfty}\le\gamma$.
\end{thm}
In Theorem~\ref{th:loi}, we generalize this LMI to a general class of infinite-dimensional systems - replacing the matrices $A,B,C,D$ with operators $\mcl A,\mcl B,\mcl C,\mcl D$ and the positive matrix variable $P$ with an operator variable $\mcl P$.

Recall that the ODE~\eqref{eq:ODE} is passive if for any input $w\in\ltwo$, we have $y\in \ltwo$ and $\ip{w}{y}_{\ltwo}\ge0$.
For ODEs, an LMI test for passivity can be formulated as follows.
\begin{thm}
If there exists a positive definite matrix $P$ such that
\begin{align}\label{eq:lmipass}
\begin{bmatrix}
A^TP+PA & PB-C^T\\
B^TP -C & -(D+D^T)\\
\end{bmatrix}\le 0
\end{align}
then for any $w\in\ltwo$ and $y\in\ltwo$ which satisfy \eqref{eq:ODE} for some $x$, $\ip{w}{y}_{\ltwo}\ge0$.
\end{thm}
In Theorem~\ref{th:loi}, we likewise generalize this LMI to infinite-dimensional systems. Having posed operator-valued feasibility tests, we next using matrices to parameterize a set of positive operators using the $PQRS$ framework and enforce positivity of such operators using LMI - See Theorem~\ref{th:positivity1}. Next, we use our new state-space framework to reduce the operator feasibility test, as applied to the PDE system in~\eqref{eq:PDE0}, to a positivity constraint on an operator of the $PQRS$ format. This feasibility test can then be verifies using LMIs, as in Theorem~\ref{th:SOS}. Numerical testing indicates the resulting bounds are not conservative in any significant sense.

%%%%%%%%%%%%%%%%%%%%%%%%%%%%%%%%%%%%%%%%%%%%%%%%%%%%%%%%%%%%%%%%%%%%%%%%%%%%%%%%%%%%%%%%%%
\section{Notation}\label{sec:notation}
We use $\mathbb{S}^m\subset \rl^{m\times m}$ to denote the symmetric matrices. We define the space of square integrable $\rl^n$-valued functions on $X$ as $\ltwo^n(X)$. $\ltwo^n(X)$ is equipped with the inner product $\langle x,y \rangle_{\ltwo} = \myint x(s)^T y(s) ds$. We also use the notation $\ip{x}{y}_{\rl}=x^Ty$ for inner product between \rl-space elements. The Sobolov space, $W^{q,n}(X):=\{x\in \ltwo^n(X) \mid \frac{\partial^k x}{\partial s^k}\in\ltwo^n(X) \text{ for all }k\le q \}$. We define the indicator function as
\mc{I}(\th) = $\left\lbrace
                            \begin{matrix} 1 \quad \th \ge 0\\
                                           0 \quad \th<0
                            \end{matrix}
                \right.
              $.
For an inner product space $X$, operator $\mc{P}:X\to X$ is called positive, if for all $x\in X$, we have $\ip{x}{\mc{P}x}_X\ge 0$. We use $\mc{P}\ge0$ to indicate that $\mcl P$ is a positive operator. We say that $\mcl P:X\to X$ is coercive if there exists some $\epsilon>0$ such that $\ip{x}{\opPa x}_X\ge \epsilon \Vert x \Vert_X^2$ for all $x \in X$.
%%%%%%%%%%%%%%%%%%%%%%%%%%%%%%%%%%%%%%%%%%%%%%%%%%%%%%%%%%%%%%%%%%%%%%%%%%%%%%%%%%%%%%%%%%
\section{LOI analogue of the Bounded-Real and Positive-Real Lemmas}\label{sec:LOIs}
Consider the abstract form of a Distributed Parameter System (DPS),
\begin{align}\label{eq:PDE}
  \dot{x}(t) &= \mathcal{A}x(t) + \mathcal{B}w(t), \nonumber\\
  y(t) &= \mathcal{C}x(t) + \mathcal{D}w(t),\qquad \qquad  x(0)=0,
\end{align}
 where, $x(t)\in X$ is the state, $y(t)\in \rl^q$ is the output and $w(t)\in\rl^m$ is the exogenous input to the system.
$\mathcal{A}:X\to Z, ~\mathcal{B}:\rl^m\to Z, ~\mathcal{C}:X\to\rl^q$ and $\mathcal{D}:\rl^m\to\rl^q$ are linear operators.

In this section, we present the conditions for passivity and \ltwo ~gain of the system \eqref{eq:PDE}.
\begin{thm} \label{th:KYP}
 %the following statements are valid.
\renewcommand{\theenumi}{\Alph{enumi}}
\begin{enumerate}
\item Suppose there exists a coercive, self-adjoint linear operator $\opPa:Z \to Z$ such that
\begin{align}\label{eq:loihinf}
  &\ip{z}{\opPa\mathcal{A}z}_Z+\ip{\mathcal{A}z}{\opPa z}_Z+\ip{z}{\opPa\mathcal{B}u}_Z\nonumber\\
  &~+\ip{\mathcal{B}u}{\opPa z}_Z\le \gamma^2\ip{u}{u}_{\rl}-\ip{\mathcal{C}z}{\mathcal{C}z}_{\rl}-\ip{\mathcal{C}z}{\mathcal{D}u}_{\rl}\nonumber\\
  &~~-\ip{\mathcal{D}u}{\mathcal{C}z}_{\rl}-\ip{\mathcal{D}u}{\mathcal{D}u}_{\rl}
\end{align} for all $z\in X\subseteq Z$ and $u\in\rl^m$. Then for any $w \in \ltwo^m([0,\infty))$ and $y \in \ltwo^q([0,\infty))$ which satisfy \eqref{eq:PDE} for some $x$, $\Vert y\Vert_{\ltwo} \le \gamma\Vert w\Vert_{\ltwo}$.
\item Suppose there exists a coercive, self-adjoint linear operator $\opPa:Z \to Z$ such that
\begin{align}\label{eq:loipass}
&\ip{z}{\opPa\mathcal{A}z}_Z+\ip{\mathcal{A}z}{\opPa z}_Z+\ip{z}{\opPa\mathcal{B}u}_Z+\ip{\mathcal{B}u}{\opPa z}_Z\nonumber\\
  &~\le\ip{\mc{C}z}{u}_{\rl}+\ip{u}{\mc{C}z}_{\rl}+\ip{\mc{D}u}{u}_{\rl}+\ip{u}{\mc{D}u}_{\rl}
\end{align} for all $z\in X\subseteq Z$ and $u\in\rl^q$. Then for any $w \in \ltwo^q([0,\infty))$ and $y \in \ltwo^q([0,\infty))$ which satisfy \eqref{eq:PDE} for some $\mbf x$, $\ip{w}{y}_{L_2}\ge 0$.
\end{enumerate}
\end{thm}
\begin{proof}
Define $V(t) = \ip{x(t)}{\opPa x(t)}_Z$. Since \mc{P} is coercive, $V(t)\ge0$.
If $x(t)$ is a solution to \eqref{eq:PDE} and $w(t)\in\rl^m$, then
\begin{align*}
  \dot{V}(t) &= \ip{x(t)}{\opPa \dot{x}(t)}_Z + \ip{\dot{x}(t)}{\opPa x(t)}_Z\\
  &= \ip{x(t)}{\opPa (\mathcal{A}x(t) + \mathcal{B}w(t))}_Z \\
  &\quad+ \ip{(\mathcal{A}x(t) + \mathcal{B}w(t))}{\opPa x(t)}_Z\\
  &= \ip{x(t)}{\opPa\mathcal{A}x(t)}_Z + \ip{x(t)}{\opPa\mathcal{B}w(t)}_Z\\
  &\quad+ \ip{\mathcal{A}x(t)}{\opPa x(t)}_Z + \ip{\mathcal{B}w(t)}{\opPa x(t)}_Z.
\end{align*}

A) Now, Inequality \eqref{eq:loihinf} implies
\begin{align*}
  \dot{V}(t) &= \ip{x(t)}{\opPa\mathcal{A}x(t)}_Z + \ip{x(t)}{\opPa\mathcal{B}w(t)}_Z\\
  &\quad+ \ip{\mathcal{A}x(t)}{\opPa x(t)}_Z + \ip{\mathcal{B}w(t)}{\opPa x(t)}_Z\\
  & \le \gamma^2\ip{w(t)}{w(t)}_{\rl}-\ip{y(t)}{y(t)}_{\rl}
\end{align*}
By the integrating the above expression with respect to time from 0 to $\infty$, we get
\begin{align*}
  \int_{0}^{\infty}\left(\dot{V}(t) + \ip{y(t)}{y(t)}_{\rl} - \gamma^2\ip{w(t)}{w(t)}_{\rl}\right) dt &\le 0.
\end{align*} Then \begin{align*}
  V(\infty) -V(0) + \ip{y}{y}_{\ltwo} - \gamma^2\ip{w}{w}_{\ltwo} &\le 0.
\end{align*}

Since $V(t)\ge0$, $\lim_{t\to\infty}V(t)\ge 0$. Also recall $x(0)=0$, so $V(0)=0$. Hence,
\begin{align*}
  \Vert y\Vert^2_{\ltwo} \le \gamma^2\Vert w\Vert^2_{\ltwo}.
\end{align*}
B) Inequality \eqref{eq:loipass} implies
\begin{align*}
  \dot{V}(t) &= \ip{x(t)}{\opPa\mathcal{A}x(t)}_Z + \ip{x(t)}{\opPa\mathcal{B}w(t)}_Z\\
  &\quad+ \ip{\mathcal{A}x(t)}{\opPa x(t)}_Z + \ip{\mathcal{B}w(t)}{\opPa x(t)}_Z\\
  & \le \ip{y(t)}{w(t)}_{\rl} + \ip{w(t)}{y(t)}_{\rl}.
\end{align*}
Integrating the above expression with respect to time from 0 to $\infty$, we get
\begin{align*}
\int_{0}^{\infty}\left(\dot{V}(t) - 2 \ip{y(t)}{w(t)}_{\rl}\right) dt &\le 0.
\end{align*} Then \begin{align*}
V(\infty) -V(0) -2\ip{y}{w}_{\ltwo} &\le 0.
\end{align*}
We recall that $V(0)=0$ and $\lim_{t\to\infty}V(t)\ge 0$. Hence $\ip{y}{w}_{\ltwo}\ge0$.
\end{proof}\addtocounter{thm}{1}
\vspace{-0.3cm}
%%%%%%%%%%%%%%%%%%%%%%%%%%%%%%%%%%%%%%%%%%%%%%%%%%%%%%%%%%%%%%%%%%%%%%
\section{Coupled PDEs in the Semigroup Framework}\label{sec:pdedef}
In the previous section, we presented conditions for passivity and \ltwo ~gain of an abstract DPS. In this section, we will focus on expressing the PDEs \eqref{eq:PDE3} in the DPS framework described in the previous section - specifically, the coupled linear PDEs of the form
\begin{align}\label{eq:PDE3}
  &\dot{x}(s,t) = A_0(s)x(s,t) + A_1(s)x_s(s,t)\notag\\
  &\qquad\qquad\qquad+A_2(s)x_{ss}(s,t) + B_1(s)w(t),\notag\\
  &y(t)= C_1 z(t) + \int_{0}^L \left(C_a(s)x(s,t)+C_b(s)x_s(s,t)\right)ds,\notag\\
  &Bz(t) = 0, \qquad x(s,0)=0,\notag\\
  &z(t)=\bmat{x(0,t)& x(L,t)& x_s(0,t)& x_s(L,t)}^T.
 \end{align}
 where $x(\cdot,t)\in X$, $y(t)\in\rl^q$ and $w\in\rl^m$.

In the semigroup framework, solutions of \eqref{eq:PDE3} also define of solutions of \eqref{eq:PDE} if $Z=\ltwo^n([0,L])$,
\begin{align*}
X &:= \{x\in \sob{[0,L]}\mid \\
&\qquad B\begin{bmatrix}x(0)& x(L)& x_s(0)& x_s(L)\end{bmatrix}^T = 0 \},
\end{align*}
and the linear operators $\mc{A}:X\to Z$, $\mc{B}:\rl^m\to Z$, $\mc{C}:X\to \rl^q$ and $\mc{D}:\rl^m\to \rl^q$ are defined as
\begin{align}\label{eq:PDE2}
  &(\mc{A}x)(s) := A_0(s) x(s) + A_1(s) x_s(s) + A_2(s) x_{ss}(s),\nonumber\\
  &(\mc{B}w)(s) := B_1(s)w, \nonumber\\
  &\mc{C}x := \myint C_3(s) x_{ss}(s) ds,
  ~\mc{D}w := D_1w,
\end{align}
where
\begin{align*}
C_3(\th) &= C_1\left(\begin{bmatrix}
I & 0\\ I & L\\ 0 &I\\0 &I
\end{bmatrix}B_c(\th) + \begin{bmatrix}
0\\ (L-\th)\\ 0 \\I
\end{bmatrix}\right)\\
&\quad+\left(\myint (C_a(s)G_1(s,\th) + C_b(s)G_2(s,\th)) ds\right),\\
G_1(s,\th) &= [I ~sI]B_c(\et) + \mc{I}(s-\th)(s-\th),\\
G_2(s,\th) &= [0 ~I]B_c(\et) + \mc{I}(s-\th),\\
B_c(\et) &= -\left(B
\begin{bmatrix}I& 0 \\I &L \\0 &I \\ 0 &I\end{bmatrix}\right)^{-1}B\begin{bmatrix}0\\L-\et\\0\\ I
\end{bmatrix}.
\end{align*}

We restrict the operators \opPa ~used in Theorem \ref{th:KYP} to a class of operators $\opPmn: Z \rightarrow Z$, parameterized by $M:\rl\to\rl^{n\times n}$ and $N_1,N_2:\rl\times\rl\to\rl^{n\times n}$ as
\begin{align}\label{eq:opdef}
  (\opPmn x)(s) &:=  M(s)x(s) + \myinta{s} N_1(s,\th) x(\th)d\th\nonumber\\
  &\qquad+\myintb{s} N_2(s,\th) x(\th)d\th.
\end{align}

We will show that for a system with operators \mc{A,~B,~C} and \mc{D} as defined in \eqref{eq:PDE2} and the class of operators \opPmn, Theorem \ref{th:KYP} can be reformulated in terms of an inequality involving operators of the form $\opPb$ defined as
\begin{align}\label{eq:opdef2}
  &\left(\opPb \vmatwo{x_1}{x_2}\right)(s) :=\nonumber\\
  &\quad \vmatwo{P x_1 + \frac{1}{L}\myint Q(s)x_2(s) ds}{Q^T(s) x_1 + \myinta{s} R_1(s,\th) x_2(\th)d\th+ \myintb{s} R_2(s,\th) x_2(\th)d\th}
\end{align}
where $P$ is a matrix, $Q$, $R_1$ and $R_2$ are matrix valued polynomials of appropriate dimensions.
 %%%%%%%%%%%%%%%%%%%%%%%%%%%%%%%%%%%%%%%%%%%%%%%%%%%%%%%%%%%%%%%%%%%%%%%%%%%%
\section{Reformulation of operator inequalities}\label{sec:reformLOI}
In Theorem \ref{th:KYP}, we saw that the problem of determining passivity and bounding the \ltwo ~gain of a DPS~\eqref{eq:PDE} parameterized by \mc{A,~B,~C} and \mc{D}  can be formulated as a feasibility test for the existence of an operator \opPa ~which satisfies the inequalities stated in the theorem. Now, we show that when the linear operators \mc{A,~B,~C} and \mc{D} ~are as defined in \eqref{eq:PDE2}, if the operator \opPa ~is parameterized by matrix valued polynomials $M$, $N_1$ and $N_2$ as described in \eqref{eq:opdef}, then inequalities in Theorem \ref{th:KYP} can be reformulated as an inequality involving operator of form \opPb defined in \eqref{eq:opdef2} and there exists a linear map from $M$, $N_1$ and $N_2$ to $P$, $Q$, $R_1$ and $R_2$.
 \begin{lem}\label{lem:loi1}
 Suppose the operators $\mc{B,~C}$ and $\mc{D}$ ~are as defined in \eqref{eq:PDE2}. Then for all $x\in X$ and $w\in\rl^m$,
 \begin{align}\label{eq:lem1eq}
     &\ip{x}{\opPmn\mc{B}w}_{\ltwo}+\ip{\mc{B}w}{\opPmn x}_{\ltwo}\nonumber\\
     &\quad+\ip{\mc{C}x}{\mc{D}w}_{\rl}+\ip{\mc{D}w}{\mc{C}x}_{\rl}+\ip{\mc{D}w}{\mc{D}w}_{\rl}\nonumber\\
     &\qquad-\gamma^2\ip{w}{w}_{\rl}=\ip{\vmatwo{w}{x_{ss}}}{\mc{P}_{\{P,Q,0,0\}}\vmatwo{w}{x_{ss}}}_{\ltwo}
 \end{align}
 where
  \begin{align*}
     P &= D_1^TD_1/L -\gamma^2/L,\\
     Q(s) &= D_1^TC_3(s) + \myint V(s,\th)^T d\th, \\
     V(s,\et)&=\Big(G(\et,s)^TM(\et)+\int_{0}^{s} G(\beta,s)^TN_1(\beta,\et)d\beta\\
     &\qquad+\int_{s}^{L} G(\beta,s)^TN_2(\beta,\et)d\beta \Big)B_1(\et),\\
     G(s,\et) &= [I ~sI]B_c(\et) + \mc{I}(s-\et)(s-\et),\\
     B_c(\et) &= -\left(B
     \begin{bmatrix}I& 0 \\I &L \\0 &I \\ 0 &I\end{bmatrix}\right)^{-1}B\begin{bmatrix}0\\L-\et\\0\\ I
     \end{bmatrix}.
     \end{align*}
 \end{lem}
 \begin{proof}
 From the boundary conditions and fundamental theorem of calculus, it can be shown that
     \begin{align}\label{eq:ftc}
     &x(s) = \myint G(s,\et)x_{ss}(\et)d\et.
     \end{align}
 We will deal with each term in the left-hand side of \eqref{eq:lem1eq} separately. Firstly,
     \begin{align*}
     &\ip{\mc{D}w}{\mc{D}w}-\gamma^2\ip{w}{w}=w^TD_1^TD_1w-\gamma^2w^Tw\\
     &\quad= \myint z(s)^T \begin{bmatrix}
     D_1^TD_1/L-\gamma^2/L & 0 \\ 0 & 0
     \end{bmatrix}z(s) ds\\
     &\quad=\ip{z}{\mc{P}_{\{P,0,0,0\}}z}
     \end{align*} where $z(s)=\vmatwo{w}{x_{ss}(s)}$.
     By substituting $x(s)$ using \eqref{eq:ftc},
     \begin{align*}
         &\ip{x}{\opPmn\mc{B}w} \\
         &\quad= \myint x(s)^T M(s)B_1(s)w ds\\
         &\quad~~+ \myint\myinta{s}x(s)^TN_1(s,\th)B_1(\th)w d\th ds\\
         &\quad~~+ \myint\myintb{s}x(s)^TN_2(s,\th)B_1(\th)w d\th ds\\
		 &\quad= \myint \myint x_{ss}(s)^T G(\et,s)^T ds M(\et)B_1(\et)w d\et\\
 		 &\quad~+ \myint \myint x_{ss}(s)^T G(\et,s)^T ds\myinta{s} N_1(\et,\th)B_1(\th)w d\th d\et\\
 		 &\quad~+ \myint \myint x_{ss}(s)^T G(\et,s)^T ds\myintb{s} N_2(\et,\th)B_1(\th)w d\th d\et\\
         &\quad= \myint x_{ss}(s)^T \left(\myint V(s,\et) d\et\right)wds =\ip{z}{\mc{P}_{\{0,Q_1,0,0\}}z}.
     \end{align*}
%     where,
%     \begin{align*}
%     V(s,\et)&=\Big(G(\et,s)^TM(\et)+\int_{0}^{s} G(\beta,s)^TN_1(\beta,\et)d\beta\\
%     &\qquad+\int_{s}^{L} G(\beta,s)^TN_2(\beta,\et)d\beta \Big)B_1(\et)\\
%     \end{align*}	
%	 Then,
%     %The term $\myint V(s,\et) d\et$ can be evaluated by numerical integration. We drop $B_1(\et)$ in the following steps for convenience. \\
%     \begin{align*}
%         &\myint V(s,\et)d\et= \myint \left([I ~\et I]B_c(\et)\right)^T M(\et)B_1(\et) d\et\\
%         &  + \int_{s}^{L} (\et-s)^T M(\et) B_1(\et)d\et+\myint [I ~\beta I]^T N_3(\beta)d\beta\\&+ \int_{s}^{L}(\beta-s)^TN_4(\beta)d\beta\\
%     \end{align*}
%     where,
%     \begin{align*}
%		N_3(\beta) &= \int_{0}^{\beta}  N_1(\beta,\et) B_1(\et) d\et+ \int_{\beta}^{L} N_2(\beta,\et)B_1(\et)d\et\\
%		N_4(\beta) &= \int_{0}^{\beta} N_1(\beta,\et)B_1(\et)  d\et+ \int_{\beta}^{L}N_2(\beta,\et)B_1(\et) d\et.
%     \end{align*}
     From \eqref{eq:PDE2},
     \begin{align*}
		\ip{\mc{C}x}{\mc{D}w}=\myint x_{ss}(s)^TC_3(s)^TD_1w ds = \ip{z}{\mc{P}_{\{0,Q_2,0,0\}}z}.
     \end{align*}
     Then
     \begin{align*}
     &\ip{x}{\opPmn\mc{B}w}_{\ltwo}+\ip{\mc{B}w}{\opPmn x}_{\ltwo}\nonumber\\
     &\quad+\ip{\mc{C}x}{\mc{D}w}_{\rl}+\ip{\mc{D}w}{\mc{C}x}_{\rl}+\ip{\mc{D}w}{\mc{D}w}_{\rl}-\gamma^2\ip{w}{w}_{\rl}\nonumber\\
     &\quad=\ip{z}{\mc{P}_{\{P,0,0,0\}}z}+ \ip{z}{\mc{P}_{\{0,Q_1,0,0\}}z}+\ip{z}{\mc{P}_{\{0,Q_2,0,0\}}z} \\
     &\quad=\ip{z}{\mc{P}_{\{P,Q,0,0\}}z}_{\ltwo}.
     \end{align*}
 \end{proof}
 \vspace{-0.5cm}
\noindent\textbf{Notation:} In the following Lemma, $\mc{L}_1$, $\mc{L}_2$ and $\mc{L}_3$ are linear maps between matrix valued polynomials that satisfy Lemmas 4, 5 and 6 of \cite{peet2018}, respectively. Detailed definition of these maps can be found in the appendix. We use these maps to establish the following Lemmas.
 \begin{lem}\label{lem:loi2}
 Suppose the operators $\mc{A}$ and $\mc{C}$ ~are as defined in \eqref{eq:PDE2}. Then for all $x\in X$ and $w\in\rl^m$,
 \begin{align*}
     &\ip{x}{\opPmn\mc{A}x}_{\ltwo}+\ip{\mc{A}x}{\opPmn x}_{\ltwo}+\ip{\mc{C}x}{\mc{C}x}_{\rl} \\
     &=\ip{\vmatwo{w}{x_{ss}}}{\mc{P}_{\{0,0,R_1,R_2\}}\vmatwo{w}{x_{ss}}}_{\ltwo}
 \end{align*}
 where
  \begin{align*}
     R_1(s,\th) &= H_1(s,\th)+H_2(\th,s)^T+C_3(s)^TC_3(\th),\\
     R_2(s,\th) &= R_1(\th,s)^T,\\
  	(H_1,H_2) &= \mc{L}_1(V_0,W_{01},W_{02})+\mc{L}_2(V_1,W_{11},W_{12}),\\
  	&\qquad+\mc{L}_3(V_2,W_{21},W_{22}),
  	\end{align*}
  	the linear maps $\mc{L}_1$, $\mc{L}_2$ and $\mc{L}_3$ are defined in the appendix and
  	\begin{align*}
  	V_i(s) = M(s)A_i(s),~W_{ij}(s,\th) = N_j(s,\th)A_i(\th)\\
  	 \forall i\in\{0,1,2\}, j\in\{1,2\}.
  \end{align*}
 \end{lem}
 \begin{proof}
We use Lemma~\ref{lem:reform2} in the Appendix to express $x$ in terms of $x_{ss}$.
 \begin{align*}
 	&\ip{x}{\opPmn\mc{A}x}\\
 	&=\ip{x}{\opPmn A_0x}+\ip{x}{\opPmn A_1x_s}\\
 	&~+\ip{x}{\opPmn A_2x_{ss}}=\ip{x_{ss}}{\mc{P}_{\{0,H_1,H_2\}} x_{ss}}.
 \end{align*}
Now, from \eqref{eq:PDE2},
 \begin{align*}
     \ip{\mc{C}x}{\mc{C}x}= \myint \myint x_{ss}(s)^T C_3(s)^TC_3(\th) x_{ss}(d\th) d\th ds.
 \end{align*}
 It follows that
  \begin{align*}
	 &\ip{x}{\opPmn\mc{A}x}+\ip{\mc{A}x}{\opPmn x}+\ip{\mc{C}x}{\mc{C}x} \\
     &= \myint\myinta{s} z(s)^T\begin{bmatrix}
	 0 & 0 \\
	 0& H_1(s,\th)+H_2(\th,s)^T
	 \end{bmatrix}z(\th)d\th ds\\
	 &\quad+\myint\myintb{s} z(s)^T\begin{bmatrix}
	 0 & 0 \\
	 0& H_2(s,\th)+H_1(\th,s)^T
	 \end{bmatrix}z(\th)d\th ds\\
     &\quad+\myint\myint z(s)^T\begin{bmatrix}
	 0 & 0 \\
	 0& C_3(s)^TC_3(\th)
	 \end{bmatrix}z(\th)d\th ds\\
     &=\ip{z}{\mc{P}_{\{0,0,R_1,R_2\}}z}.
     \end{align*}
 \end{proof}
 \vspace{-0.5cm}
In the following two theorems, we combine the preceding two lemmas to reformulate the inequalities of Theorem \ref{th:KYP} for the coupled linear PDE system \eqref{eq:PDE3} defined in Section \ref{sec:pdedef} in terms of an inequality using an operator of the form \opPb.
 \addtocounter{thm}{-1}
\begin{thm}\label{th:loi}
	 Suppose there exists polynomials M: $\mathbb{R}\rightarrow \mathbb{R}^{n\times n}$ and $N_1,~N_2$: $\mathbb{R}\times \mathbb{R} \rightarrow \mathbb{R}^{n\times n}$ such that \opPmn ~is coercive and $\opPb\le0$,
% 	 $\ip{\vmatwo{w}{x}}{\opP \vmatwo{w}{x}}_{\ltwo}\leq 0$ for all  $\vmatwo{w}{x} \in \mathbb{R}^m \times \sob$
where $P$, $Q$, $R_1$ and $R_2$ are as defined in Lemmas \ref{lem:loi1} and \ref{lem:loi2}.
     For any $w \in \ltwo^m([0,\infty))$, if $x(t)\in X$ and $y(t)$ satisfy \eqref{eq:PDE3}, then $\norm{y}_{\ltwo}\le\gamma\norm{w}_{\ltwo}$.
     \end{thm}
	 \begin{proof}
	 From Lemmas \ref{lem:loi1} and \ref{lem:loi2},
     \begin{align*}
  &\ip{x}{\opPmn\mathcal{A}x}_{\ltwo}+\ip{\mathcal{A}x}{\opPmn x}_{\ltwo}\nonumber\\
  &\quad+\ip{x}{\opPmn\mathcal{B}w}_{\ltwo}+\ip{\mathcal{B}w}{\opPmn x}_{\ltwo}\nonumber\\
  &\quad~~+\ip{\mathcal{C}x}{\mathcal{C}x}_{\rl}+\ip{\mathcal{C}x}{\mathcal{D}w}_{\rl}+\ip{\mathcal{D}w}{\mathcal{C}x}_{\rl}+\ip{\mathcal{D}w}{\mathcal{D}w}_{\rl}\nonumber\\
  &\quad\quad~~-\gamma^2\ip{w}{w}_{\rl} \\
  &\quad= \ip{z}{\mc{P}_{\{P,0,0,0\}}z}+\ip{z}{\mc{P}_{\{0,0,R_1,R_2\}}z}+\ip{z}{\mc{P}_{\{0,Q,0,0\}}z}\\
  &\quad=\ip{\vmatwo{w}{x_{ss}}}{\opPb\vmatwo{w}{x_{ss}}}_{\ltwo}\le0.
     \end{align*}
     Then Inequality \eqref{eq:loihinf} is satisfied and hence $\norm{y}_{\ltwo}\le\gamma\norm{w}_{\ltwo}$.
     \end{proof}
     \vspace{-0.5cm}
\begin{thm}\label{th:loi2}
     Suppose there exists polynomials $M:\mathbb{R}\rightarrow \mathbb{R}^{n\times n}$ and $N_1,~N_2$: $\mathbb{R}\times \mathbb{R} \rightarrow \mathbb{R}^{n\times n}$, such that \opPmn ~is coercive and $\opPb\le0$,
where $P$, $Q$, $R_1$ and $R_2$ are defined as
     \begin{align} \label{eq:passive}
     P &= -(D_1+D_1^T)\nonumber,\\
     Q(s) &= -C_3(s) +\myint V(s,\th)^Td\th,\nonumber\\
     R_1(s,\th) &= H_1(s,\th)+H_2(\th,s)^T,\nonumber\\
     R_2(s,\th) &= R_1(\th,s)^T,
     \end{align}
     where $V(s,\th)$ is defined in Lemma \ref{lem:loi1},
     \begin{align*}
     (H_1,H_2) &= \mc{L}_1(V_0,W_{01},W_{02})+\mc{L}_2(V_1,W_{11},W_{12})\\
     &\qquad+\mc{L}_3(V_2,W_{21},W_{22}),
	 \end{align*}
	the linear maps $\mc{L}_1$, $\mc{L}_2$ and $\mc{L}_3$ are defined in the appendix and
	\begin{align*}
	V_i(s) = M(s)A_i(s),~W_{ij}(s,\th) = N_j(s,\th)A_i(\th)\\
	\forall i\in\{0,1,2\}, j\in\{1,2\}.
	\end{align*}
     For any $w\in\ltwo^q([0,\infty))$, if $x(t)\in X$ and $y(t)\in\rl^q$ satisfy \eqref{eq:PDE3}, then $\ip{w}{y}_{\ltwo}\ge0$.
\end{thm}
\begin{proof}
	Again, we use the results from Lemmas \ref{lem:loi1} and \ref{lem:loi2}, to express $x$ in terms of $x_{ss}$. We deal with each term separately. First,
     \begin{align*}
         &-(\ip{\mc{D}w}{w}+\ip{w}{\mc{D}w}) = -w^T(D_1^T+D_1)w\\
         &= \myint z(s)^T \begin{bmatrix}
         -(D_1^T+D_1) & 0 \\ 0 & 0
         \end{bmatrix}z(s) ds=\ip{z}{\mc{P}_{\{P,0,0,0\}}z},
     \end{align*} where z(s) = $\vmatwo{w}{x_{ss}(s)}$.\\
     Next, we have
     \begin{align*}
	 &\ip{x}{\opPmn\mc{B}w+\ip{\mc{B}w}{\opPmn x}}-\ip{\mc{C}x}{w}\\
	 &\quad-\ip{w}{\mc{C}x}\\
     &=\myint x_{ss}(s)^T\left(\myint V(s,\th) d\th\right) wds \\
     &\quad+\myint w^T\left(\myint V(\th,s)^Tds\right) x_{ss}(\th) d\th \\
     &\quad-\myint x_{ss}(s)^TC_3(s)^Tw ds-\myint w^TC_3(s)x_{ss}(s)ds\\
     &=\myint z(s)^T\begin{bmatrix}
	 0 & -C_3(s)+\myint V(s,\th)^T\\
	 *^T & 0
	 \end{bmatrix}z(s) ds\\
     &=\ip{z}{\mc{P}_{\{0,Q,0,0\}}z}.
	 \end{align*}
	 From Lemma \ref{lem:loi2},
	 \begin{align*}
     &\ip{x}{\opPmn\mc{A}x}+\ip{\mc{A}x}{\opPmn x}\\
     &=\ip{z}{\mc{P}_{\{0,0,R_1,R_2\}}z}
 \end{align*}
 where
 \begin{align*}
     R_1(s,\th) = H_1(s,\th)+H_2(\th,s)^T ~\text{and} ~R_2(s,\th) = R_1(\th,s)^T.
 \end{align*}
     Then
     \begin{align*}
         &\ip{x}{\opPmn\mc{A}x}+\ip{\mc{A}x}{\opPmn x}\\
         &~+\ip{x}{\opPmn\mc{B}w}+\ip{\mc{B}w}{\opPmn x}\\
          &~~-(\ip{\mc{C}x}{w}+\ip{w}{\mc{C}x}+\ip{\mc{D}w}{w}+\ip{w}{\mc{D}w})\\
          &~~~=\ip{z}{\mc{P}_{\{0,0,R_1,R_2\}}z}+\ip{z}{\mc{P}_{\{0,Q,0,0\}}z}\\
          &~~~~~+\ip{z}{\mc{P}_{\{P,0,0,0\}}z}\\
          &~~~=\ip{z}{\opPb z}_{\ltwo}\le0.
     \end{align*}
    We conclude that Inequality \eqref{eq:loipass} is satisfied and hence the system is passive.
\end{proof}
\noindent\textbf{Notation:} For convenience, we define two new linear maps. Specifically, if  $P$, $Q$, $R_1$, $R_2$, $M$, $N_1$ and $N_2$ satisfy Theorem \ref{th:loi} then we say
\begin{align*}
     \{P,Q,R_1,R_2\} = \mc{L}_4(M,N_1,N_2).
\end{align*}
Likewise, if  $P$, $Q$, $R_1$, $R_2$, $M$, $N_1$ and $N_2$ satisfy Theorem \ref{th:loi2} then we say
\begin{align*}
     \{P,Q,R_1,R_2\} = \mc{L}_5(M,N_1,N_2).
\end{align*}
%%%%%%%%%%%%%%%%%%%%%%%%%%%%%%%%%%%%%%%%%%%%%%%%%%%%%%%%%%%%%%%%%%%%%%%%%%%%%%%%%%%%%%%%%%%%%%%%%%
\section{Enforcing positivity of operators of form \opPb}\label{sec:opIneq}
In Theorem~\ref{th:KYP}, we showed that the problem of determining passivity and $H_\infty$ gain of an abstract Distributed Parameter System (DPS) - parameterized by the operators $\mcl A, \mcl B, \mcl C$ and $\mcl D$ - could be formulated as a convex feasibility problem of the existence of operator $\mcl P$ which satisfies certain operator inequalities. In Equation~\eqref{eq:PDE2}, we showed that coupled PDE systems with inputs and outputs could be cast in a DPS framework by defining the operators $\mcl A, \mcl B, \mcl C$ and $\mcl D$ for this class of systems. In Equation~\eqref{eq:opdef}, we used matrix-valued functions $M,N_1$ and $N_2$ to parameterize a class of operators, denoted $\mcl P_{\{M,N_1,N_2\}}$ acting on the state space defined by the system of coupled PDEs. Next, in Theorem~\ref{th:loi} and \ref{th:loi2}, we showed that, using these definitions and parameterization of variables, the feasibility conditions of Theorem~\ref{th:KYP} could be expressed as positivity of $\mcl P_{\{M,N_1,N_2\}}$ and negativity of an operator $\opPb$ parameterized by $P$, $Q$, $R_1$ and $R_2$ as defined in Equation~\eqref{eq:opdef2} where if $M$, $N_1$ and $N_2$ are polynomials, there is a linear map from the coefficients in $M$, $N_1$ and $N_2$ to the elements of $P$ and the coefficients of the polynomials $Q$, $R_1$, and $R_2$. In the following two theorems, we show how to use LMI constraints to enforce positivity of the operators $\mcl P_{\{M,N_1,N_2\}}$ and $\opPb$, respectively. These results will be used in Theorem~\ref{th:SOS} to give an SDP representation of Theorem~\ref{th:KYP} as applied to the coupled PDE system in~\eqref{eq:PDE2}.
\begin{thm}\label{th:positivity1}
	For any functions $Z_1:X\to\rl^{d_1\times n}$, $Z:X\times X\rightarrow\mathbb{R}^{d_2\times n}$, suppose there exists a matrix $T\ge0$ such that
	\setlength{\abovedisplayskip}{\mylen}
	\setlength{\belowdisplayskip}{\mylen}
	\begin{align}
	M &= Z_1(s)^TT_{11}Z_1(s),\nonumber\\
	N_1(s,\th) &=Z_1(s)^TT_{12}Z(s,\th)+ Z(\th,s)^TT_{32}Z(\th)\\
	&\hspace{-1.25cm}+\myintb{s}Z(\beta,s)^TT_{22}Z(\beta,\th)d\beta+\int_{\th}^{s}Z(\beta,s)^TT_{32}Z(\beta,\th)d\beta\nonumber\\
	&+\myinta{\th}Z(\beta,s)^TT_{33}Z(\beta,\th)d\beta,\nonumber\\
	N_2(s,\th) &=N_1(\th,s)^T,
	\end{align}
	where
	\begin{align*}
	T= \begin{bmatrix}
	T_{11} & T_{12} & T_{13} \\
	T_{21} & T_{22} & T_{23} \\
	T_{31} & T_{32} & T_{33} \\
	\end{bmatrix}.
	\end{align*}\\
	Then for the operator \opPmn ~as defined in \eqref{eq:opdef}, $\opPmn\ge0$.\\
\end{thm}
\vspace{-0.5cm}
See~\cite{peet2014LMI} for a proof.
\begin{thm}\label{th:positivity}
    For any function $Z(s,\th):X\times X\rightarrow\mathbb{R}^{d_2\times n}$, suppose there exists a matrix %$T\in \mathbb{S}^{(m+2n)} $ such that
    $T\ge0$ such that
    \setlength{\abovedisplayskip}{\mylen}
    \setlength{\belowdisplayskip}{\mylen}
   \begin{align}
    P &= T_{11},\nonumber\\
    Q(\th) &= \myintb{\th} T_{12}Z(s,\th)ds + \myinta{\th} T_{13}Z(s,\th)ds, \nonumber\\
    R_1(s,\th) &=\myintb{s}Z(\beta,s)^TT_{22}Z(\beta,\th)d\beta\nonumber\\
    &\hspace{-1.2cm}+\int_{\th}^{s}Z(\beta,s)^TT_{32}Z(\beta,\th)d\beta+\myinta{\th}Z(\beta,s)^TT_{33}Z(\beta,\th)d\beta,\nonumber\\
    R_2(s,\th) &=\myintb{\th}Z(\beta,s)^TT_{22}Z(\beta,\th)d\beta\nonumber\\
    &\hspace{-1.2cm}+\int_{s}^{\th}Z(\beta,s)^TT_{23}Z(\beta,\th)d\beta+\myinta{s}Z(\beta,s)^TT_{33}Z(\beta,\th)d\beta,\nonumber\\
    \label{eq:TH}
   \end{align}
   where
   \begin{align*}
   T= \begin{bmatrix}
   T_{11} & T_{12} & T_{13} \\
   T_{21} & T_{22} & T_{23} \\
   T_{31} & T_{32} & T_{33} \\
   \end{bmatrix}.
   \end{align*}\\
%   is partitioned such that $T_{11}\in\mathbb{S}^m,~ T_{22},~T_{33}\in\mathbb{S}^n$.
   Then for the operator \opPb ~as defined in \eqref{eq:opdef}, $\opPb\ge0$.
   \end{thm}
   \begin{proof}
   \setlength{\abovedisplayskip}{\mylen}
    \setlength{\belowdisplayskip}{\mylen}
   Since $T\ge0$, we can define a square root of T as $U = \begin{bmatrix} U_1&U_2&U_3\end{bmatrix}$.
   \setlength{\arraycolsep}{4pt}
   \begin{align*}\label{eq:Teq}
   T &= \begin{bmatrix}
   T_{11} & T_{12} & T_{13}\\
   *^T & T_{22} & T_{23} \\
   *^T & *^T & T_{33} \\
   \end{bmatrix} = U^T U=\begin{bmatrix}
   U_1^TU_1 & U_1^TU_2 & U_1^TU_3 \\
   *^T & U_2^TU_2 & U_2^TU_3 \\
   *^T & *^T & U_3^TU_3\\
   \end{bmatrix}.
   \end{align*}
   Let define $v(s) = U_1x_1 + (\Psi x_2)(s)$, where
\begin{align*}
(\Psi x_2)(s)&=\myinta{s} U_2 Z(s,\th) x_2(\th) d\th+\myintb{s} U_3 Z(s,\th) x_2(\th) d\th.
\end{align*}
 Then
   \begin{align*}
       &\ip{v}{v}_{\ltwo}= \ip{U_1x_1}{U_1x_1}_{\ltwo}+\ip{U_1x_1}{(\Psi x_2)}_{\ltwo}\\
   &~~~+\ip{(\Psi x_2)}{U_1x_1}_{\ltwo}+\ip{(\Psi x_2)}{(\Psi x_2)}_{\ltwo}\\
   &\quad=\ip{\vmatwo{x_1}{x_2}}{\opPb \vmatwo{x_1}{x_2}}_{\ltwo}\ge 0.
   \end{align*}
%
%   From the theorem statement,
%   \begin{align*}
%   &\ip{U_1x_1}{U_1x_1}_{\ltwo}=\myint (U_1x_1)^T (U_1x_1) ds= \myint x_1^TT_{11}x_1ds\\
%   &\quad =\myint x_1^T P x_1 ds
%   \end{align*}
%   and
%   \begin{align*}
%   &\ip{U_1x_1}{(\Psi x_2)(s)}_{\ltwo}=\myint (U_1x_1)^T (\Psi x_2)(s) ds\\
%   &= \myint x_1^T \left(\begin{matrix}
%   \myintb{s} U_1^T U_2Z(\th,s)d\th\\
%   \quad+ \myinta{s} U_1^T U_3Z(\th,s)d\th
%   \end{matrix}\right)x_2(s)ds\\
%   &=\myint x_1^T \left(\begin{matrix}
%   \myintb{s} T_{12}Z(\th,s)d\th\\ +\myinta{s} T_{13}Z(\th,s)d\th
%   \end{matrix}\right)x_2(s)ds\\
%   &=\myint x_1^TQ(s)x_2(s)ds.
%   \end{align*}
%Using Lemma 7 of \cite{peet2014LMI}, it can be shown that
%\begin{align*}
%& \ip{(\Psi x_2)(s)}{(\Psi x_2)(s)}_{\ltwo}=\myint (\Psi x_2)(s)^T (\Psi x_2)(s) ds \\
%&= \myint\myinta{s} x_2^T(s) R_1(s,\th)x_2(\th)d\th ds \\
%   &~~+ \myint\myintb{s} x_2^T(s) R_2(s,\th)x_2(\th)d\th ds.
%\end{align*}
%Then,
%\begin{align*}
%   &\ip{v(s)}{v(s)}_{\ltwo} \\
%   &\quad= \myint \vmatwo{x_1}{x_2(s)}^T\begin{bmatrix}
%   P & Q(s) \\ Q^T(s) &0
%   \end{bmatrix}\vmatwo{x_1}{x_2(s)} ds \\
%   &\quad~+ \myint\myinta{s} x_2^T(s) R_1(s,\th)x_2(\th)d\th ds \\
%   &\quad~~+ \myint\myintb{s} x_2^T(s) R_2(s,\th)x_2(\th)d\th ds\\
%   &\quad=\ip{\vmatwo{x_1}{x_2}}{\opPb \vmatwo{x_1}{x_2}}_{\ltwo}.
%   \end{align*}
Hence $\opPb \ge 0$.
\end{proof}
For convenience, we define the following two sets.
\begin{align*}
\Phi_1 &:= \{\{P,Q,R_1,R_2\} : P,~Q,~R_1 ~\text{and}~R_2 \\
&\quad\text{satisfy the conditions of Theorem \ref{th:positivity}}\}.
\end{align*}
\begin{align*}
\Phi_2 &:= \{\{M,N_1,N_2\} \mid M, N_1 ~\text{and} ~N_2 \\
&\quad\text{satisfy the conditions of Theorem \ref{th:positivity1}}\}.
\end{align*}
%%%%%%%%%%%%%%%%%%%%%%%%%%%%%%%%%%%%%%%%%%%%%%%%%%%%%%%%%%%%%%%%%%%%%%%%%
\section{An SOS formulation for \hinfty ~analysis}\label{sec:Hinf}
In this section, we consolidate Lemmas and Theorems from Sections \ref{sec:reformLOI} and \ref{sec:opIneq} to arrive at the LMI equations that are sufficient to test passivity and find the bound on \ltwo ~gain of the PDE \eqref{eq:PDE3}.
\begin{thm}\label{th:SOS}
Suppose there exists $\epsilon>0$, $\gamma>0$, matrix-valued polynomials $M:\rl\to\rl^{n\times n}$, and $N_1,N_2:\rl\times\rl\to\rl^{n\times n}$ such that
\begin{align*}
(M-\epsilon I, N_1, N_2) &\in \Phi_2.\end{align*}
Then for all $x(t)\in X$, $y\in \ltwo^q([0,\infty))$ and $w\in \ltwo^m([0,\infty))$ which satisfy \eqref{eq:PDE3},
\begin{enumerate}
    \item if
    \begin{align*}
\{P,Q,R_1,R_2\} =\mc{L}_4(M,N_1,N_2)
\end{align*}
such that \begin{align*}
\{-P,-Q,-R_1,-R_2\} &\in \Phi_1,
\end{align*} then $\Vert y\Vert_{\ltwo} \le \gamma\Vert w\Vert_{\ltwo}$.
\item if $m=q$ and% there exists $P\in\mathbb{S}^m$, $Q:\rl\to\rl^{m\times n}$ and $R_1,R_2:\rl\times\rl\to\rl^{n\times n}$ such that,
\begin{align*}
\{P,Q,R_1,R_2\} =\mc{L}_5(M,N_1,N_2)
\end{align*} such that \begin{align*}
\{-P,-Q,-R_1,-R_2\} &\in \Phi_1,
\end{align*} then $\ip{y}{w}_{\ltwo}\ge0$
\end{enumerate}
\end{thm}
\begin{proof}
Suppose that $V(x)=\ip{x}{\opPmn x}_{\ltwo}\ge\epsilon\Vert x\Vert_{\ltwo}^2$.

1) Since $\{P,Q,R_1,R_2\} =\mc{L}_4(M,N_1,N_2)$, Theorem \ref{th:loi} is satisfied.
Then
\begin{align*}
  &\ip{x}{\opPmn\mathcal{A}x}_{\ltwo}+\ip{\mathcal{A}x}{\opPmn x}_{\ltwo}\nonumber\\
  &\quad+\ip{x}{\opPmn\mathcal{B}w}_{\ltwo}+\ip{\mathcal{B}w}{\opPmn x}_{\ltwo}\nonumber\\
  &\quad\quad+\ip{\mathcal{C}x}{\mathcal{C}x}_{\rl}+\ip{\mathcal{C}x}{\mathcal{D}w}_{\rl}+\ip{\mathcal{D}w}{\mathcal{C}x}_{\rl}+\ip{\mathcal{D}w}{\mathcal{D}w}_{\rl}\nonumber\\
  &\quad\quad\quad-\gamma^2\ip{w}{w}_{\rl} = \ip{\vmatwo{w}{x_{ss}}}{\opPb \vmatwo{w}{x_{ss}}}_{\ltwo}.
\end{align*} \normalsize
From the conditions of the theorem, $\opPb\le 0$. Consequently, Inequality \eqref{eq:loihinf} from part (A) of Theorem \ref{th:KYP} is satisfied and hence $\Vert y\Vert_{\ltwo} \le \gamma\Vert w\Vert_{\ltwo}$.

2) The proof for the second part of this theorem is quite similar to the first.
Since $\{P,Q,R_1,R_2\} =\mc{L}_5(M,N_1,N_2)$, Theorem \ref{th:loi2} is satisfied.
From Theorem \ref{th:loi2},
\begin{align*}
&\ip{x}{\opPmn\mathcal{A}x}_{\ltwo}+\ip{\mathcal{A}x}{\opPmn x}_{\ltwo}\nonumber\\
  &~+\ip{x}{\opPmn\mathcal{B}w}_{\ltwo}+\ip{\mathcal{B}w}{\opPmn z}_{\ltwo}\\&~~-(\ip{\mc{C}x}{w}_{\rl}+\ip{x}{\mc{C}w}_{\rl}+\ip{\mc{D}w}{w}_{\rl}+\ip{w}{\mc{D}w}_{\rl})\\
  &~~~= \ip{\vmatwo{w}{x_{ss}}}{\opPb \vmatwo{w}{x_{ss}}}_{\ltwo}.
\end{align*}
Since $\opPb\le0$, Inequality \eqref{eq:loipass} from part (B) of Theorem \ref{th:KYP} is satisfied and hence $\ip{w}{y}_{\ltwo}\ge0$.
\end{proof}
%The equations obtained from Theorem \ref{th:SOS}, when solved while minimizing $\gamma$, provides an estimate for \hinfty ~bound for the system.

%%%%%%%%%%%%%%%%%%%%%%%%%%%%%%%%%%%%%%%%%%%%%%%%%%%%%%%%%%%%%%%%%%%%%%%%%%%%%%%%%%%%%%%%%%%%%%%%
\section{Numerical Simulations and Validation}\label{sec:numeri}
Algorithm presented in Theorem~\ref{th:SOS} was implemented in MATLAB. We compare the estimate of $H_\infty$ norm bound obtained by using numerical discretization with the estimate from our method, for several PDE systems. In all cases, referring to~\eqref{eq:PDE2} we use the following values
\begin{align*}
	B_1(s)=1, ~C_1=0, ~C_a(s)=1, ~C_b(s)=0 ~\text{and} ~D_1=0.
\end{align*}
\subsection{Example 1}\label{ex:spat}
Consider the system shown below. In \cite{peet2018}, it was shown to be stable for $\lambda< 4.65$.
\begin{align*}
u_t(s,t)  = A_0(s) u(s,t) + &A_1(s)u_s(s,t) \\&+ A_2(s)u_{ss}(s,t)+w(t)\\
u(0,t) = 0\quad &\quad u_s(L,t) = 0
\end{align*}
\begin{align*}
A_0(s) &= (-0.5s^3+1.3s^2-1.5s+0.7+\lambda)\\
A_1(s) &= (3s^2-2s), \quad A_2(s) = (s^3-s^2+2)
\end{align*}

\setlength{\belowcaptionskip}{0pt}
\setlength{\abovecaptionskip}{0pt}

\begin{figure}[!h]
    \centering
    \begin{subfigure}[b]{0.5\textwidth}
        \centering
        \includegraphics[width=\textwidth]{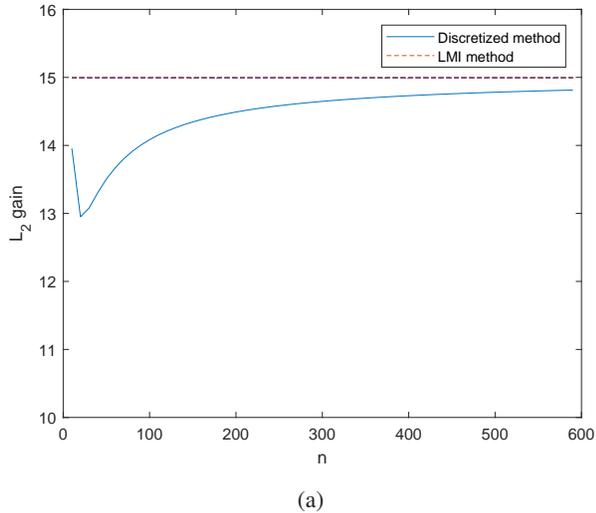}
        \caption{}
        \label{fig:gainpoints}
    \end{subfigure}
    \begin{subfigure}[b]{0.5\textwidth}
        \centering
        \includegraphics[width=\textwidth]{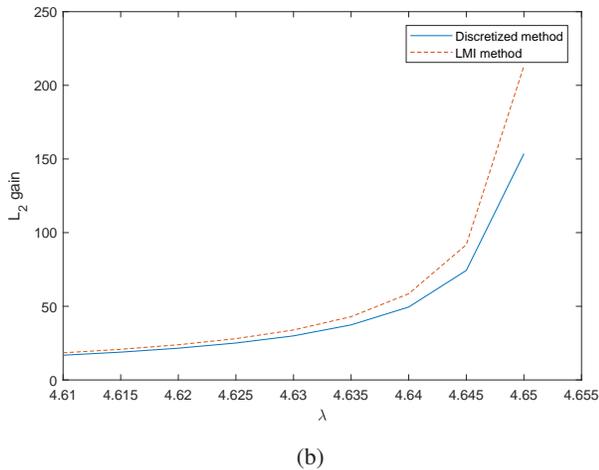}
        \caption{}
        \label{fig:gainparam}
    \end{subfigure}
    \caption{For the PDE system in Example \ref{ex:spat}: (a) Mesh size vs \ltwo ~gain obtained from Theorem~\ref{th:SOS} and spatial discretization, (b)  value of the parameter, $\lambda$ vs \ltwo ~gain obtained from Theorem~\ref{th:SOS} and spatial discretization}
\end{figure}

 Fig. \ref{fig:gainpoints} shows the variation of an estimate of the \ltwo ~gain obtained from spatial discretization while varying mesh size. At at mesh size of 600, we had an \ltwo ~gain of 14.82 (LMI bound was 14.99). Although this example obtained the largest residual gap of all examples at 3\%, this residual is likely due to our naive method of discretization and not conservatism in Theorem~\ref{th:SOS}. Fig. \ref{fig:gainparam} shows the bounds obtained when the system parameter $\lambda$ is varied. Using higher degree polynomials shows minor change in the \ltwo ~gain bound, typically of the order $10^{-6}$. This suggests that relatively low degree polynomials give tight bounds.
\subsection{}\label{ex:3eq}
For the PDE systems listed below, we compare the \ltwo ~gain bounds obtained by our algorithm and finite-difference dicretization method in Table I.
\begin{enumerate}[label=\textit{B}.\arabic*:]
	\item Following PDE is stable for $\lambda\leq\pi^2$. \begin{align*}
		&u_t(s,t)  = \lambda u(s,t) + u_{ss}(s,t)\\
		&u(0,t)=0,\quad u(L,t)=0.
	\end{align*}
	\item Following PDE is stable for $\lambda\leq2.467$.\begin{align*}
	&u_t(s,t)  = \lambda u(s,t) + u_{ss}(s,t)\\
	&u(0,t)=0,\quad u_s(L,t)=0.
	\end{align*}
	\item The following coupled PDE was shown to be stable for $R<21$ in \cite{ahmadi2016} .
	\begin{align*}
	u_t(s,t)  = \begin{bmatrix}
	0 & 0 & 0\\ s & 0 & 0\\ s^2 &-s^3 &0
	\end{bmatrix}& u(s,t) + \frac{1}{R}u_{ss}(s,t)+w(t)\\
	u(0,t) = 0\quad &\quad u(L,t) = 0
	\end{align*}
\end{enumerate}
\begin{table}[h]\label{tab:list}
	\begin{center}
	\begin{tabular}{|c|c|c|c|}
		\hline
		          &LMI method   &Discretized method &Parameter\\\hline
		B.1&8.214&8.253 &$\lambda=0.98\pi^2$\\\hline
		B.2&12.03&12.31 &$\lambda=2.4$\\\hline
		B.3&3.9738&3.9708 &$R=20$\\\hline
	\end{tabular}
	\caption{A bound on \ltwo ~gain using different methods.}
	\end{center}
\end{table}
%\subsection{Example 2}\label{ex:3eq}
% The values for \hinfty ~bound from SOS method and numerical discretization methods are 3.974 and 3.971 respectively. Computation time to solve the optimization problem was 5.6s and degree of the polynomials used in parametric form of Lyapunov function was 2.
\subsection{Example 3} \label{ex:diss}
Consider,
\begin{align*}
u_{t,i}(s,t) = \lambda u_i(s,t)+&\sum_{k=1}^{i}u_{ss,k}(s,t)+w(t)\\
u(0,t) = 0\quad &\quad u(L,t) = 0.
\end{align*}
This example was tailored to test the time complexity of the algorithm proposed. We use the value $\lambda=0.5\pi^2$ for all $i$. CPU time of the algorithm for different number of coupled PDEs is tabulated in Table II.
\begin{table}[h]\label{tab:heateq}
%\begin{center}
\begin{tabular}{|c|c|c|c|c|c|c|c|}
\hline
i          &1   &2   &3   &4   &5   &10 &20 \\\hline
CPU time(s)&0.60&1.45&5.22&13.7&36.5&2317 &27560 \\\hline
\end{tabular}
\caption{Runtime for the system of equations for increasing number of coupled PDEs, i. Refer Example \ref{ex:diss}}
%\end{center}
\end{table}
%%%%%%%%%%%%%%%%%%%%%%%%%%%%%%%%%%%%%%%%%%%%%%%%%%%%%%%%%%%%%%%%%%%%%%%%%%%%%%%%
\vspace{-5pt}
\section{CONCLUSIONS}\label{sec:conclu}
In this paper, we proposed a method to prove passivity and obtain bounds for the \ltwo -gain of coupled linear PDEs with domain distributed disturbances using the LMI framework. The method presented does not use discretization. The bounds and properties obtained are prima facie provable. The numerical results indicate there is little, if any conservatism in the result.

\section{Acknowledgments}
This work was supported by National Science Foundation under grants No. 1739990, 1538374 and by Office of Naval Research Award N00014-17-1-2117.
\bibliographystyle{IEEEtran}
\bibliography{References}

% \vfill

\appendix
\section{}
\addtocounter{thm}{1}
We restate main result from Lemmas 4, 5 and 6 from \cite{peet2018}. These results are used in the proof of Lemma \ref{lem:loi2} in Section \ref{sec:reformLOI}.

\begin{defn}For given matrix-valued functions $M$, $N_1$ and $N_2$ and given matrix $B \in \R^{2n \times 4n}$ of row rank $2n$, we say that
	\begin{align*}
	(R_1,R_2)&=\mcl L_1(M,N_1,N_2),\quad	(Q_1,Q_2)=\mcl L_2(M,N_1,N_2)\\
	(T_1,T_2)&=\mcl L_3(M,N_1,N_2)
	\end{align*}
	if {\small
		\begin{align*}
		R_1(s,\theta)&=E_1(s,\theta)+E_3(s,\theta),&
		R_2(s,\theta)&=E_2(s,\theta) +E_3(s,\theta),\\
		Q_1(s,\theta)&=F_1(s,\theta)+F_3(s,\theta),&
		Q_2(s,\theta)&=F_2(s,\theta) +F_3(s,\theta),\\
		T_1(s,\theta)&=G_1(s,\theta)+G_3(s,\theta),&
		T_2(s,\theta)&=G_2(s,\theta) +G_3(s,\theta),
		\end{align*}}%
\noindent where  {\small
		\begin{align*}
		&E_1(s,\theta)=\int_s^b (\eta-s)  N_1(\eta,\theta)d\eta\\
		&E_2(s,\theta)=(\theta-s)M(\theta)+\int_\theta^b (\eta-s) N_1(\eta,\theta)d\eta\\
		&\hspace{4.5cm} +\int_s^\theta(\eta-s) N_2(\eta,\theta)d\eta \\
		&E_3(s,\theta) =Y_1(s,\theta)\\
%		\end{align*}
%		\\\vspace{-20pt}
%		\begin{align*}
		&F_1(s,\theta)= \int_s^b\left(  (\eta-s)F_4(\theta,\eta)+ F_5(s,\eta) \right)d \eta\\
		&F_2(s,\theta)=\int_\theta^b \left( (\eta-s) F_4(\theta,\eta) + F_5(s,\eta)\right)d \eta\\
		&F_3(s,\eta) =\int_a^b B_a(\zeta,s)^T Y_2(\zeta) B_b(\eta) d \zeta +\int_\eta^b Y_1(s,\zeta)d \zeta\\
		&\hspace{4cm} +\int_s^b (\zeta-s)Y_2(\zeta) d \zeta B_b(\eta)\\
		&F_4(\theta,\eta)= M(\eta)  +  \int_\theta^\eta   N_1(\eta,\zeta)   d \zeta\\
		&F_5(s,\eta)= \int_s^\eta  (\zeta - s)  N_2(\zeta,\eta)    d \zeta\\
		&G_1(s,\theta)=\int_s^b \left((\eta-s) G_4(\theta,\eta) +G_5(s,\theta,\eta) \right) d \eta \\
		&G_2(s,\theta)=\int_\theta^b \left((\eta-s) G_4(\theta,\eta)+ G_5(s,\theta,\eta) \right) d \eta\\
		&G_3(s,\theta)=\int_a^b B_a(\eta,s)^T Y_3(\eta,\theta)d\eta\\
		&\hspace{1.3cm}+\int_\theta^b (\eta-\theta) Y_1(s,\eta) d \eta+\int_s^b (\eta-s) Y_3(\eta,\theta) d\eta\\
		&G_4(\theta,\eta)= (\eta-\theta) M(\eta)  +  \int_\theta^\eta  (\zeta-\theta) N_1(\eta,\zeta)  d \zeta \\
		&G_5(s,\theta,\eta)=\int_s^\eta  (\zeta-s) (\eta-\theta)  N_2(\zeta,\eta)   d \zeta\\
%		\end{align*}}%
%	{\small
%		\begin{align*}
		&Y_1(s,\eta)=B_a(\eta,s)^T M(\eta) + \int_\eta^b  B_a(\theta,s)^T N_1(\theta,\eta)d\theta\\
		&\hspace{2.7cm}+\int_a^\eta B_a(\theta,s)^T N_2(\theta,\eta)d\theta\\
		&Y_2(\zeta)=M(\zeta)  + \int_a^\zeta    N_1(\zeta,\theta)   d \theta + \int_\zeta^b  N_2(\zeta,\theta)   d \theta\\
		&Y_3(\zeta,\eta)=  M(\zeta)B_a(\zeta,\eta) + \int_a^\zeta   N_1(\zeta,\theta)  B_a(\theta,\eta)  d \theta \\
		&\hspace{2.7cm}+  \int_\zeta^b N_2(\zeta,\theta) B_a(\theta,\eta)  d \theta\\
%		\end{align*}}%
%	where
%	\begin{align*}
	&B_a(s,\eta)=B_4(s)(b-\eta)+B_5(s),\\
	&B_b(\eta)=B_6(b-\eta) +  B_7\\
	&\bmat{B_6 & B_7}  =\bmat{0& I} B_3,\\
	& \bmat{B_4(s) & B_5(s)}  =\bmat{I&(s-a)I} B_3\\
	&B_3=B_2^{-1}B \bmat{0 & 0\\ I & 0\\ 0 & 0\\ 0& I},\qquad B_2=B\bmat{I &0\\I&(b-a)I\\0&I\\0 &I}.
	\end{align*}}
\end{defn}

\begin{lem}\label{lem:reform2}For given matrix-valued functions $M$, $N_1$ and $N_2$ and given matrix $B \in \R^{2n \times 4n}$ of row rank $2n$, suppose that $(R_1,R_2)=\mcl L_1(M,N_1,N_2)$, $(Q_1,Q_2)=\mcl L_2(M,N_1,N_2)$ and $(T_1,T_2)=\mcl L_3(M,N_1,N_2)$. Then for any $\mbf x \in X$ where $\mbf x $ is as defined in Eqn.~\eqref{eq:PDE2}, we have that
	\begin{align*}	\ip{\mbf x}{\mcl P_{\{M,N_1,N_2\}}\mbf x_{ss}}&= \ip{\mbf x_{ss}}{\mcl P_{\{0,R_1,R_2\}}\mbf x_{ss}}\vspace{-2mm}\\
	\ip{\mbf x}{\mcl P_{\{M,N_1,N_2\}}\mbf x_s}&= \ip{\mbf x_{ss}}{\mcl P_{\{0,Q_1,Q_2\}}\mbf x_{ss}}\vspace{-2mm}\\
	\ip{\mbf x}{\mcl P_{\{M,N_1,N_2\}}\mbf x}&= \ip{\mbf x_{ss}}{\mcl P_{\{0,T_1,T_2\}}\mbf x_{ss}}\vspace{-2mm}
	\end{align*}
\end{lem}

%%%%%%%%%%%%%%%%%%%%%%%%%%%%%%%%%%%%%%%%%%%%%%%%%%%%%%%%%%%%%%%%%%%%%%%%%%%%%%%%
\end{document}